\documentclass [a4paper,10pt]{article}
\usepackage [utf8] {inputenc}
\usepackage {amsmath}
\usepackage {amsthm}
\usepackage {amssymb}
\usepackage {amsfonts}
\newtheorem* {theorem}{Theorem}
\title {Primes in arithmetic progressions on average II}
\author {Tomos Parry}
\date {}
\begin {document}
\maketitle
\begin {abstract}
A deep conjecture of Montgomery and Soundararajan on the distribution of prime numbers in short intervals of length $h$ says that the third moment is bounded by $\ll h^{\frac {3}{2}-c}$ for some $c>0$. There is in the literature some conditional evidence towards this conjecture whilst in the first article to this series we gave the first instance of unconditional evidence in the form of a bound corresponding to $\ll h^{7/5+o(1)}$. In this article we push the exponent down to $\ll h^{1+o(1)}$ which more or less is expected to be best possible.
\end {abstract}
\hspace {1mm}
\\ 
\begin {center}
\section {- \hspace {3mm}Introduction}
\hspace {1mm}
\end {center}
\hspace {1mm}
\\ Let
\[ N^\epsilon \leq h\leq N^{1-\epsilon }\]
and define
\[ \Lambda (n):=\left \{ \begin {array}{ll}\log p&\text { if $n$ is a power of $p$}\\ 0&\text { if not}\end {array}\right \} \hspace {10mm}\psi (t):=\sum _{t<n\leq t+h}\Lambda (n).\] 
A conjecture of Montgomery and Soundararajan \cite {montsound} says that $\psi (t)$ is normally distributed with mean $h$ and variance $h\log (N/h)$, i.e. that
\begin {equation}\label {conjecture}
\frac {1}{N}\int _1^N\Big (\psi (t)-h\Big )^kdt\hspace {1mm}\sim \hspace {1mm}\mu _k\left (h\log (N/h)\right )^{k/2}
\end {equation}
where $\mu _k$ are the moments of the standard normal distribution. For $k=2$ this is expected to be equivalent to the strong pair correlation conjecture. These conjectures are far out of reach but may be accessible conditionally - indeed Montgomery and Soundararajan prove \eqref {conjecture} on the assumption of the Riemann hypothesis and the Hardy-Littlewood $k$-tuple conjecture (see their Theorem 3).
\\
\\ Let $R_k(h)$ be as in (8) of \cite {montsound}. According to (19) in \cite {montsound} and the sentences following it, we would expect the LHS of \eqref {conjecture} to amount to $NR_k(h)$ (on the above conjectures). These $R_k(h)$ have been studied further by Bloom and Kuperberg \cite {bloomkuperberg} so that the current state of affairs has 
(see Theorems 2 of both \cite {montsound} and \cite {bloomkuperberg})
\[ R_{k}(h)=\mu _k\left (-h\log h+Ah\right )^{k/2}+\mathcal O\left (h^{\frac {k-1}{2}+\epsilon }\right )\]
which of course for odd $k$ reads
\[ R_{k}(h)\ll h^{\frac {k-1}{2}+\epsilon }.\]
For $(a,q)=1$ let
\begin {equation}\label {primes}
E_x(q,a):=\sum _{p\leq x\atop {p\equiv a(q)}}\log p-\frac {x}{\phi (q)}\hspace {20mm}\phi (q)\hspace {1mm}:=\hspace {1mm}\sum _{a=1\atop {(a,q)=1}}^q1.
\end {equation}
A ``$q$-analogue" of \eqref {conjecture} for odd moments would read
\[ \frac {1}{q}\sum _{a=1\atop {(a,q)=1}}^qE_x(q,a)^k\ll \left (\frac {x}{q}\right )^{\frac {k-1}{2}+\epsilon }\]
and a (again highly conditional) result of this kind has been established by Sun-Kai Leung \cite {leung}. Note that $q$-analogues are typically no easier than their short-interval versions.
\\
\\ These results are important because they provide conditional evidence for \eqref {conjecture}. However, as far as we are aware, unconditional evidence is sparse. A result of Hooley \cite {hooley} amounts to
\[ \frac {1}{Q^2}\sum _{q\leq Q}\sum _{a=1\atop {(a,q)=1}}^qE_x(q,a)^3\ll \left (\frac {x}{Q}\right )^{3/2}e^{-c\sqrt {\log (x/Q)+2}}+\frac {x^2}{(\log x)^A}\hspace {10mm}\text {for }Q\leq \frac {x}{\log x}\] 
which detects $o(h^{3/2})$ for the $q$-analogue but which is still only a dent into the above discussion's predicted 
\[ \ll \left (\frac {x}{Q}\right )^{1+\epsilon }.\]
In \cite {tomos}, however, we showed a power-saving unconditionally, i.e. an exponent $3/2-\delta $, and in this article we push the exponent down to what is essentially best possible. We prove (unconditionally)
\begin {theorem} 
For $Q\leq x/(\log x)^2$ and any $A>0$
\[ \sum _{q\leq Q}\phi (q)\sum _{a=1\atop {(a,q)=1}}^qE_x(q,a)^3\ll _AQ^3\left (\frac {x}{Q}\right )^{1+\epsilon }+\frac {x^3}{(\log x)^A}.\]
\end {theorem}
{}\hspace {1mm}
\\ Should we want to study the fourth moment with a dispersion method (this is in fact our main motivation) then, as indicated in \cite {tomos}, this theorem offers the best input.
\\
\\ Now to the proof. When the letters $l,l'$ are in context write $\tilde l:=l+l'$. If the letters $f,g$ are multiplicative functions (specifically only these letters) and if $\Delta \in \mathbb N$ define
\begin {eqnarray}\label {bound2}
f_\Delta (n)&:=&\prod _{p^N||n\atop {p\nmid \Delta }}f(p^N)\hspace {6mm}\text { so that }\hspace {6mm}f_\Delta (dn)\hspace {1mm}=\hspace {1mm}f_\Delta (d)f_{d\Delta }(n)
\end {eqnarray}
and similarly for $g_\Delta $. By $\int _{(c)}$ we mean a contour integral running from $c-i\infty $ to $c+i\infty $. Let
\begin {eqnarray}\label {s}
r&:=&\frac {1}{p}\left \{ \begin {array}{ll}1&\text { if }p>2\\ 0&\text { if }p=2\end {array}\right \} +\mathcal O\left (\frac {1}{p^2}\right )\hspace {10mm}f_\Delta (n):=\prod _{p|n\atop {p\nmid \Delta }}\left (1+r\right )\hspace {10mm}f_\Delta =:g_\Delta \star 1\notag 
\\ \mathcal Q_\Delta (s)&:=&\prod _{p\nmid \Delta }\left (1+\frac {2r}{p}+\frac {rY}{p}\left (p-1+r\right )\right )\notag 
\\ &=&\prod _{p\nmid \Delta }\frac {1}{1-Y/p}\prod _{p\nmid \Delta }\left (1-\frac {Y}{p}\right )\left (\underbrace {1+\frac {2r}{p}+\frac {rY}{p}\left (p-1+r\right )}_{1+Y\left (r-\frac {1}{p}\right )+\mathcal O\left (\frac {1+Y^2}{p^2}\right )}\right )\hspace {10mm}Y\hspace {1mm}:=\hspace {1mm}p^{-s}\notag 
\\ \mathcal E_\Delta (X)\hspace {1mm}&:=&\hspace {1mm}2\int _{(-1/4)}\frac {\mathcal Q_\Delta (s)X^{s+4}ds}{s(s+3)(s+4)}\hspace {1mm}\ll \hspace {1mm}X^4+1\notag 
\\ \mathcal S_\Delta (X)&:=&\sum _{q}g_\Delta (q)\sum _{l+l'\leq X\atop {\tilde l\equiv 0(q)\atop {(ll',q)=1}}}(X-\tilde l)^2\tilde lf_\Delta (l)f_\Delta (l')
\end {eqnarray}
and let $\mathcal M_\Delta $ denote any quantity of the form $\alpha _\Delta X^5+X^4\left (\beta _\Delta \log X+\gamma _\Delta \right )$ with $\alpha _\Delta ,\beta _\Delta ,\gamma _\Delta \ll 1$; note
\begin {eqnarray}
g_\Delta (q)=\left \{ \begin {array}{ll}g_1(q)&\text { if }(q,\Delta )=1\\ 0&\text { if }(q,\Delta )>1\text { or if $q$ is $\square$-free}\end {array}\right \} \hspace {5mm}\text { and }\hspace {2mm}g_1(p)=r\hspace {5mm}\text { so }\hspace {5mm}g_\Delta (q)\ll q^{\epsilon -1}.\label {bound}
\end {eqnarray}
The paper is divided into sections 2,3,4. In \S 2 we will prove that our theorem follows from
\begin {equation}\label {seren}
\mathcal S_\Delta (X)\hspace {1mm}=\hspace {1mm}\mathcal M_\Delta (X)+\mathcal E_\Delta (X)+\mathcal O_\epsilon \left (X^{3+\epsilon }\right ).
\end {equation}
Implicit in \cite {hooley} is \eqref {seren} with an error like $\mathcal O\left (X^{7/2}e^{-\sqrt {\log X}}\right )$ and with a different main term. As explained in \cite {tomos}, the form of Hooley's main term means that the error term can't be improved without RH. In \cite {tomos} we sorted out the main term and consequently proved \eqref {seren} with error $\mathcal O(X^{17/5+\epsilon })$. For the much stronger error term in this article, we investigate $\mathcal S_\Delta $ more closely by taking into account cancellation in $\chi $ in \eqref {l} below - previously this was bounded absolutely in both \cite {hooley} and \cite {tomos}. Let $R_\Delta (n)$ be multiplicative with
\begin {eqnarray}\label {r}
R_\Delta (p^\alpha )&:=&\left \{ \begin {array}{ll}r-1/p&\text { if }\alpha =1\\ -r/p&\text { if }\alpha =2\\ 0&\text { if }\alpha >2\end {array}\right \} \hspace {10mm}p\nmid \Delta \notag 
\\ &:=&\left \{ \begin {array}{ll}-1/p&\text { if }\alpha =1\\ 0&\text { if }\alpha >1\end {array}\right \} \hspace {10mm}p|\Delta .
\end {eqnarray}
For a Dirichlet character $\chi $ mod $(q)$ let
\begin {eqnarray}
\mathcal F_\chi (s)&:=&\sum _{l=1}^\infty \frac {f_\Delta (l)\chi (l)}{l^s}\label {ff}
\\ &=&L_\chi (s)\prod _{p\nmid \Delta }\left (1+r\chi (p)Y\right )\hspace {10mm}Y\hspace {1mm}:=\hspace {1mm}p^{-s}\label {hsv}%
\\ &=&L_\chi (s)L_\chi (s+1)\sum _{n=1}^\infty \frac {R_\Delta (n)}{n^s}\hspace {10mm}\text {with}\hspace {10mm}R_\Delta (n)\ll n^{\epsilon -1}\label {extension}
\end {eqnarray}
and let\footnote {that $\mathcal S_\Delta ^>=0$ may not have escaped the beady-eyed, but bear with us!}
\begin {align}
\mathcal L_\Delta (q)&:=\sum _{l+l'\leq X\atop {\tilde l\equiv 0(q)\atop {(ll',q)=1}}}(X-\tilde l)^2\tilde lf_\Delta (l)f_\Delta (l')\label {ll}
\\ \mathcal S_\Delta ^<&:=\sum _{q\leq X}g_\Delta (q)\mathcal L_\Delta (q)&\mathcal S_\Delta ^>\hspace {1mm}&:=\sum _{q>X}\mathcal L_\Delta (q)\label {ss}
\\ w_s(p)&:=1+\frac {2r}{p}+\frac {Yr^2}{p}&
\mathcal P_\Delta (s)&:=\prod _{p\nmid \Delta }w_s(p)
\hspace {10mm}\theta _s(q):=\frac {\phi (q)}{q}\prod _{p|q}\frac {1}{w_s(p)}
\label {p}
\\ \mathcal A_\Delta ^<(X)&:=\sum _{q\leq X}\frac {g_\Delta (q)\mathcal F_{\chi _0}(1)^2}{\phi (q)}&\mathcal A_\Delta ^>(X)&:=\hspace {1mm}\mathcal P_\Delta (1)\sum _{q>X}\frac {g_\Delta (q)\theta _1(q)}{q}\label {a}
\\ \mathcal Q_\Delta ^<(s)&:=\sum _{q=1}^\infty \frac {g_\Delta (q)}{\phi (q)}\mathcal F_{\chi _0}(1)\mathcal F_{\chi _0}(s)&\mathcal Q_\Delta ^>(s)&:=\zeta (s)\mathcal P_\Delta (s)\sum _{q=1}^\infty \frac {g_\Delta (q)\theta _{s}(q)}{q^s}\label {q}
\\ \mathcal E_\Delta ^<&:=2\int _{(-1/2)}\frac {\mathcal Q_\Delta ^<(s)X^{s+4}ds}{s(s+2)(s+3)(s+4)}&\mathcal E_\Delta ^>&:=2\int _{(1/2)}\frac {\mathcal Q_\Delta ^>(s)X^{s+4}ds}{(s+2)(s+3)(s+4)}\label {e}
\\ R_s&:=\sum _{q\leq X}\frac {g_\Delta (q)\theta _{s}(q)}{q^s}&\hspace {10mm}\mathcal R_\Delta &:=2\int _{(-1)}\pi ^{s-1/2}\frac {\Gamma \left (\frac {1}{2}-\frac {s}{2}\right )}{\Gamma \left (\frac {s}{2}\right )}\frac {\zeta (1-s)\mathcal P_\Delta (s)R_sX^{s+4}ds}{(s+2)(s+3)(s+4)}.&
\end {align}
Using \eqref {hsv} we have
\begin {eqnarray*}
\mathcal F_{\chi _0}(s)&=&\zeta (s)\prod _{p\nmid \Delta }\left (1+rY\right )\prod _{p|q}\frac {1}{1+rY}\prod _{p|q}(1-Y)\hspace {1mm}:=\hspace {1mm}\zeta (s)\mathcal F_\Delta ^*(s)\psi _s(q)\Delta _s(q)\hspace {10mm}\text {for }(q,\Delta )=1
\end {eqnarray*}
so from \eqref {bound} a quick check shows 
\begin {eqnarray*}
\mathcal Q_\Delta ^<(s)&=&\zeta (s)\mathcal F_\Delta ^*(s)\prod _{p\nmid \Delta }\left (1+\frac {r\psi _1(p)\Delta _1(p)\psi _s(p)\Delta _s(p)}{\phi (p)}\right )\hspace {1mm}=\hspace {1mm}\mathcal Q_\Delta (s)
\\ \mathcal Q_\Delta ^>(s)&=&\zeta (s)\mathcal P_\Delta (s)\prod _{p\nmid \Delta }\left (1+rY\theta _{s}(p)\right )\hspace {1mm}=\hspace {1mm}\mathcal Q_\Delta (s)
\end {eqnarray*}
so that $2\mathcal E_\Delta ^<+\mathcal E_\Delta ^>=\mathcal E_\Delta $. Also \eqref {hsv} and \eqref {p} say
\begin {eqnarray*}
\frac {\mathcal F_{\chi _0}(1)^2}{\phi (q)}=\frac {\mathcal P_\Delta (1)\theta _1(q)}{q}\hspace {10mm}\text {for }(q,\Delta )=1
\end {eqnarray*}
so from \eqref {a} and \eqref {bound} we see that $\mathcal A_\Delta ^<(X)+\mathcal A_\Delta ^>(X)$ is independent of $X$. We will prove in \S 3
\begin {eqnarray}\label {1}
\mathcal S_\Delta ^<\hspace {1mm}=X^5\mathcal A^<(X)+X^4\left (\beta \log X+\gamma \right )+2\mathcal E_\Delta ^<+\mathcal R_\Delta +\mathcal O\left (X^{3+\epsilon }\right )\hspace {10mm}\text {for some $\beta ,\gamma \ll 1$}
\end {eqnarray}
and in \S 4 
\begin {equation}\label {2}
\mathcal S_\Delta ^>\hspace {1mm}=\hspace {1mm}X^5\mathcal A^>(X)+\mathcal E_\Delta ^>-\mathcal R_\Delta +\mathcal O\left (X^{3+\epsilon }\right )
\end {equation}
which together imply \eqref {seren}.
\begin {center}
\section {- \hspace {3mm}Proof that \eqref {seren} $\implies $ our theorem}{\hspace {1mm}}
\end {center}
Recall the definitions of $f_\Delta ,g_\Delta ,\mathcal S_\Delta (X)$ in \eqref {s} so that in particular
\begin {eqnarray}
\mathcal J_\Delta ^*(X)&:=&\sum _{l+l'\leq X}(X-\tilde l)^2\tilde lf_\Delta (l)f_\Delta (l')f_\Delta (\tilde l)\notag  
\\ &=&\sum _{q}g_\Delta (q)\sum _{l+l'\leq X\atop {\tilde l\equiv 0(q)}}(X-\tilde l)^2\tilde lf_\Delta (l)f_\Delta (l')\label {split}
\\ &=&
\sum _{d\leq X}d^3g_\Delta (d)f_\Delta (d)^2\mathcal S_{d\Delta }(X/d)\label {calculus}
\end {eqnarray}
after splitting the inner sum of \eqref {split} according to the value of $d:=(l,q)=(l',q)$ and using \eqref {bound2}. With $r$ as in \eqref {s} define 
\begin {eqnarray}\label {largesieve}
I(p)&:=&\left \{ \begin {array}{ll}-r(1+3r+r^2)&\text { if $p$ is odd}\\ 2&\text { if $p=2$ }\end {array}\right \} \hspace {10mm}I(\Delta )\hspace {1mm}:=\hspace {1mm}\prod _{p|\Delta }\left \{ \begin {array}{ll}I(p)&\text { if $\Delta $ is $\square $-free}\\ 0&\text { if not }\end {array}\right \} \ll \Delta ^{\epsilon -1}\notag 
\\ \mathcal J^*(X)&:=&\sum _{\Delta }\Delta ^3I(\Delta )\mathcal J_\Delta ^*(X/\Delta ).
\end {eqnarray}
Recall the definition of $\mathcal M_\Delta $ in the sentence containing \eqref {bound}. Assuming \eqref {seren} we know $\mathcal S_{d\Delta }(X/d\Delta )$ is, for some $\alpha _{d\Delta },\beta _{d\Delta },\gamma _{d\Delta }\ll 1$,
\begin {eqnarray*}
\alpha _{d\Delta }\left (\frac {X}{d\Delta }\right )^5+\left (\frac {X}{d\Delta }\right )^4\left (\beta _{d\Delta }\log \left (\frac {X}{d\Delta }\right )+\gamma _{d\Delta }\right )+\mathcal E_{d\Delta }\left (\frac {X}{d\Delta }\right )+\mathcal O\left (\left (\frac {X}{d\Delta }\right )^{3+\epsilon }\right )
\end {eqnarray*}
so \eqref {calculus}, \eqref {largesieve} and the bounds \eqref {bound}, \eqref {s} give, for some $\alpha ,\beta ,\gamma $,
\begin {eqnarray}\label {nos}
\mathcal J^*(X)&=&X^5\sum _{d\Delta \leq X}\frac {I(\Delta )g_\Delta (d)f_\Delta (d)^2\alpha _{d\Delta }}{(d\Delta )^2}+...\notag 
\\ &&+\hspace {4mm}\sum _{d\Delta \leq X}I(\Delta )(d\Delta )^3g_\Delta (d)f_\Delta (d)^2\mathcal E_{d\Delta }\left (\frac {X}{d\Delta }\right )+\mathcal O\left (X^3\sum _{d\Delta \leq X}\left |I(\Delta )g_\Delta (d)f_\Delta (d)^2\right |\right )\notag 
\\ &=&\alpha X^5+\beta X^4\log X+\gamma X^4+\mathcal E(X)+\mathcal O\left (X^{3+\epsilon }\right )
\end {eqnarray} 
with
\begin {eqnarray*}
\mathcal E&:=&2\int _{(-1/2)}\frac {\zeta (s)}{s(s+3)(s+4)}\left (\sum _{d,\Delta =1}^\infty \frac {I(\Delta )g_\Delta (d)f_
\Delta (d)^2}{(d\Delta )^{s+1}}\mathcal Q_{d\Delta }(s)\right )X^{s+4}ds\hspace {10mm}(\text {$\mathcal Q_{\Delta }(s)$ as in \eqref {s}})
\end {eqnarray*}
so from \eqref {bound} and the definitions \eqref {s}, \eqref {largesieve} this $d,\Delta $-series is
\begin {eqnarray*}
&=&\left (\prod _{p}q_s(p)\right )\sum _{d,\Delta =1\atop {(d,\Delta )=1}}^\infty \frac {I(\Delta )g_1(d)f_1(d)^2}{(d\Delta )^{s+1}q_s(d\Delta )}\hspace {10mm}q_s(p)\hspace {1mm}:=\hspace {1mm}1+\frac {2r}{p}+\frac {rY}{p}(p-1+r)
\\ &=&\prod _p\left (q_s(p)+\frac {Y}{p}\Big (I(p)+g_1(p)f_1(p)^2\Big )\right )
\\ &=&\prod _p\left (1+\frac {Y}{p-1}\right ).
\end {eqnarray*}
But Lemma 1 of \cite {tomos} forces the main term in \eqref {nos} to be the same as in that lemma, so that in fact Lemma 1 holds with error $\mathcal O(X^{3+\epsilon })$. We can\footnote {a reader uncomfortable with taking the word of a preprint can always refer to \cite {hooley}, which has the corresponding statement but with error $\mathcal O\left (X^{7/2}e^{-c\sqrt {\log X}}\right )$, but there will be some work in collecting together various formulas in the paper, perhaps starting from (112)} then follow the rest of the proof in \cite {tomos} word for word but with $7/5$ replaced by $1$.
\begin {center}
\section {- \hspace {3mm}Proof of \eqref {1}}
\end {center}
By character orthogonality, the usual Perron type representations and \eqref {ff}
\begin {eqnarray}\label {l}
&&\sum _{q\leq X}g_\Delta (q)\sum _{l+l'\leq X\atop {\tilde l\equiv 0(q)\atop {(ll',q)=1}}}(X-\tilde l)^2lf_\Delta (l)f_\Delta (l')\notag 
\\ &&\hspace {10mm}=\hspace {4mm}
\Gamma (3)\sum _{q\leq X}\frac {g_\Delta (q)}{\phi (q)}\sum _{\chi (q)}\chi (-1)\underbrace {\int _{(1)}\int _{(1)}\frac {\mathcal F_{\chi }(s)\mathcal F_{\overline \chi }(w)\Gamma (s)\Gamma (w+1)X^{s+w+3}dwds}{\Gamma (s+w+4)}}_{=:\hspace {1mm}\mathcal C_\chi }.
\end {eqnarray} 
Let 
$\epsilon >0$. By \eqref {hsv} we can move these integrals to the left to get 
\begin {eqnarray}\label {banana}
\mathcal C_{\chi }&=&\underbrace {\frac {\mathcal F_{\chi }(1)^2X^5}{\Gamma (6)}}_{\chi =\chi _0}+\underbrace {\mathcal F_{\chi }(1)}_{\chi =\chi _0}\int _{(\epsilon )}\mathcal F_{\chi }(s)\overbrace {\frac {\Gamma (s)+\Gamma (s+1)}{\Gamma (s+5)}}^{\hspace {1mm}=\hspace {1mm}\frac {1}{s(s+2)(s+3)(s+4)}\hspace {1mm}=:\hspace {1mm}g(s)}X^{s+4}ds
\notag 
\\ &&+\hspace {4mm}
\int _{(\epsilon )}\int _{(\epsilon )}\frac {\mathcal F_{\chi }(s)\mathcal F_{\chi }(w)\Gamma (s)\Gamma (w+1)X^{s+w+3}dwds}{\Gamma (s+w+4)}.
\end {eqnarray}
The first two terms are (when present)
\begin {eqnarray}\label {be}
&&\frac {\mathcal F_{\chi }(1)^2X^5}{\Gamma (6)}+\mathcal F_{\chi }(1)\left (Res_{s=0}\left \{ \mathcal F_{\chi }(s)g(s)X^{s}\right \} +\int _{(-1/4)}\mathcal F_{\chi }(s)g(s)X^{s}ds\right )X^4=:R_q+\int I_q\notag 
\end {eqnarray}
so from \eqref {ss}, \eqref {ll} and \eqref {l}
\begin {eqnarray}
\mathcal S_\Delta ^<&=&2\Gamma (3)\sum _{q\leq X}\frac {g_\Delta (q)}{\phi (q)}\left (R_q+\int I_q+\sum _{\chi (q)}\chi (-1)\int _{(\epsilon )}\int _{(\epsilon )}\frac {\mathcal F_{\chi }(s)\mathcal F_{\chi }(w)\Gamma (s)\Gamma (w+1)X^{s+w+3}dwds}{\Gamma (s+w+4)}\right )\notag 
\\ &=:&2\Gamma (3)\sum _{q\leq X}\frac {g_\Delta (q)}{\phi (q)}\left (R_q+\int I_q+\int _{(\epsilon )}\int _{(\epsilon )}\frac {f_{s,w}(q)\Gamma (s)\Gamma (w+1)X^{s+w+3}dwds}{\Gamma (s+w+4)}\right )\label {ll}
\\ &=:&2\mathcal R+2\Gamma (3)\int \! \! \! \int \mathcal D.\label {dwr}
\end {eqnarray}
Recall the definitions of $\mathcal A^<(X),\mathcal E_\Delta ^<$ from \eqref {a}, \eqref {e} and the sentence containing \eqref {s}. The definitions of $\mathcal F_\chi (s)$ and $g(s)$ (in \eqref {banana} and \eqref {ff}) imply that 
\[ R_q=X^5\mathcal F_{\chi _0}(1)^2+X^4\left (\beta _q\log X+\gamma _q\right )\hspace {10mm}\text {some }\beta _q,\gamma _q\]
whilst from \eqref {extension} we get comfortably $I_q\ll X^4$ so
\begin {eqnarray}\label {qq}
\mathcal R&=&=\frac {2\Gamma (3)}{\Gamma (6)}\sum _{q\leq X}\frac {g_\Delta (q)}{\phi (q)}+\Gamma (3)\sum _{q=1}^\infty \frac {g_\Delta (q)}{\phi (q)}\left (X^4\left (\beta _q\log X+\gamma _q\right )+\int I_q\right )+\mathcal O\left (X^4\sum _{q>X}\frac {1}{q\phi (q)}\right )\hspace {10mm}\text {(using \eqref {bound})}\notag 
\\ &=&X^5\mathcal A^<(X)+X^4\left (\beta \log X+\gamma _q\right )+\Gamma (3)\mathcal M_\Delta +\mathcal E_\Delta ^<+\mathcal O\left (X^{3+\epsilon }\right )
\end {eqnarray}
and now we turn to $\mathcal D$. Until \eqref {eira} assume that $s,w$ have real part $\epsilon $. Let
\begin {equation}\label {ll}
\mathcal L_{\chi }(s,w):=L_\chi (s)L_\chi (s+1)L_{\overline {\chi }}(w)L_{\overline {\chi }}(w+1),
\end {equation}
let $\epsilon (\chi )$ be the usual root number in the functional equation for $L_\chi (s)$, and let $\chi (-1)=:(-1)^k$. By moving\footnote {the $\Gamma \left (\frac {w+s+k}{2}\right )$ factor means the integrand is rapidly decaying} the integral
\[ \int _{(2)}\frac {\Gamma \left (\frac {w+s+k}{2}\right )}{\Gamma \left (\frac {s+k}{2}\right )}L_\chi (w+s)\frac {dw}{wX^w}\]
to the line $\sigma =-2$, using the functional equation, changing variables from $-w\mapsto w$, and integrating the series for $L_{\chi }$ termwise, it may be established that for primitive $\chi \not =\chi _0$
\begin {eqnarray}\label {hofrenydd}
L_{\chi }(s)&=&\sum _{n=1}^\infty \frac {\chi (n)}{n^s}W_s\left (\frac {\sqrt d}{nX}\right )+\left (\frac {d}{\pi }\right )^{1/2-s}\epsilon (\chi )g_{\chi }(s)\sum _{n=1}^\infty \frac {\overline {\chi }(n)}{n^{1-s}}W_{1-s}\left (\frac {\sqrt dX}{n}\right )
\end {eqnarray}
where for $X>0$
\begin {eqnarray*}
W_u(X)&:=&\int _{(c)}\frac {\Gamma \left (\frac {w+u+k}{2}\right )}{\Gamma \left (\frac {u+k}{2}\right )}\frac {X^wdw}{w}\ll _cX^c\hspace {10mm}c>0
\\ g_\chi (s)&:=&\frac {\Gamma \left (\frac {1-s+k}{2}\right )}{\Gamma \left (\frac {s+k}{2}\right )}.
\end {eqnarray*}
Let $u$ have real part $>0$. By pushing the integral to $-\mathfrak R\mathfrak e(u)+\epsilon $ we have
\begin {eqnarray}\label {diagonal}
W_u(X)&=&1+\mathcal O\left (X^{-\mathfrak R\mathfrak e(u)+\epsilon }\right )
\end {eqnarray}
whilst by choosing a very large $c$ 
\[ \sum _{n=1}^\infty \left |W_u\left (\frac {1}{nd^\epsilon }\right )\right |
\ll \frac {1}{d^{100}}\hspace {10mm}\text {and}\hspace {10mm}\sum _{n>d^{1+2\epsilon }}W_u\left (\frac {d^{1+\epsilon }}{n}\right )
\ll \frac {1}{d^{100}}\]
so \eqref {hofrenydd} becomes
\begin {eqnarray*}
L_{\chi }(s)&=&\left (\frac {d}{\pi }\right )^{1/2-s}\epsilon (\chi )g_{\chi }(s)\sum _{n\leq d^{1+2\epsilon }}\frac {\overline {\chi }(n)}{n^{1-s}}W_{1-s}\left (\frac {d^{1+\epsilon }}{n}\right )+\mathcal O\left (\frac {1}{d^{99}}\right )
\end {eqnarray*}
and similarly
\begin {eqnarray*}
L_{\chi }(s+1)&=&\sum _{n\leq d^{1+2\epsilon }}\frac {\chi (n)}{n^{s+1}}W_{s+1}\left (\frac {d^{1+\epsilon }}{n}\right )+\mathcal O\left (\frac {1}{d^{99}}\right )
\end {eqnarray*}
so altogether \eqref {ll} becomes (with the obvious interpretation for ``$\mathbf n\leq d$")
\begin {eqnarray}\label {hwyr}
\mathcal L_{\chi }(s,w)&=&\left (\frac {d}{\pi }\right )^{1-s-w}\underbrace  {\epsilon (\chi )\epsilon (\overline {\chi })}_{=\hspace {1mm}1}\underbrace  {g_{\chi }(s)g_{\overline {\chi }}(w)}_{=:\hspace {1mm}\mathcal G_{\chi }(s,w)}\sum _{\mathbf n\leq d^{1+2\epsilon }}\frac {\chi \left (-\overline {n_1}n_2n_3\overline {n_4}\right )}{n_1^{1-s}n_2^{s+1}n_3^{1-w}n_4^{w+1}}\notag 
\\ &&\hspace {1mm}\times \hspace {3mm}\underbrace {W_{1-s}\left (\frac {d^{1+\epsilon }}{n_1}\right )W_{s+1}\left (\frac {d^{1+\epsilon }}{n_2}\right )W_{1-w}\left (\frac {d^{1+\epsilon }}{n_3}\right )W_{w+1}\left (\frac {d^{1+\epsilon }}{n_4}\right )}_{=:\mathcal W_{s,w}(\mathbf n)}\hspace {2mm}+\hspace {2mm}\mathcal O\left (\frac {1}{d^{50}}\right )\hspace {10mm}
\end {eqnarray}
giving with \eqref {extension}
\begin {eqnarray}\label {n}
\mathcal F_\chi (s)\mathcal F_{\overline {\chi }}(w)&=&\left (\frac {d}{\pi }\right )^{1-s-w}\mathcal G_{\chi }(s,w)\sum _{\mathbf n\leq d^{1+2\epsilon }\atop {(\mathbf n,d)=1}}\frac {R_\Delta (n_5)R_\Delta (n_6)}{n_5^sn_6^w}\frac {\chi \left (-\overline {n_1}n_2n_3\overline {n_4}n_5\overline {n_6}\right )}{n_1^{1-s}n_2^{s+1}n_3^{1-w}n_4^{w+1}}\mathcal W_{s,w}(\mathbf n)+\mathcal O\left (\frac {1}{d^{50}}\right )\notag 
\\ &=:&\left (\frac {d}{\pi }\right )^{1-s-w}\mathcal G_{\chi }(s,w)\sum _{\mathbf n\atop {(\mathbf n,d)=1}}\mathcal N_{s,w}(\mathbf n)\chi (-N)+\mathcal O\left (\frac {1}{d^{50}}\right ).
\end {eqnarray}
A character $\chi $ mod $q$ of conductor $d$ has $\chi =\chi _0\chi ^*$ for a primitive character $\chi ^*$ mod $d$ and the principal character $\chi _0$ mod $q$; in particular 
\[ L_\chi (s)=L_{\chi ^*}(s)\sum _{l|q}\frac {\mu (l)\chi ^*(l)}{l^s}.\]
Then \eqref {ll}, the last sentence and \eqref {n} imply
\begin {eqnarray}
f_{s,w}(q)&=&\sum _{d|q}\sideset {}{^*}\sum _{\chi ^*(d)}\mathcal F_{\chi _0\chi ^*}(s)\mathcal F_{\overline {\chi _0\chi ^*}}(w)\notag 
\\ &=&\sum _{d|q}\left (\frac {d}{\pi }\right )^{1-s-w}\sideset {}{^*}\sum _{\chi ^*(d)}\mathcal G_{\chi ^*}(s,w)\sum _{\mathbf l|q\atop {(\mathbf l,d)=1}}\frac {\mu (\mathbf l)}{l_1^sl_2^{s+1}l_3^{w}l_4^{w+1}}\sum _{\mathbf n\atop {(\mathbf n,d)=1\atop {(n_5n_6,q)=1}}}\mathcal N_{s,w}(\mathbf n)\chi ^*\left (-Nl_1l_2\overline {l_3}\overline {l_4}\right )\notag 
\\ &=:&\sum _{d|q}\left (\frac {d}{\pi }\right )^{1-s-w}\sideset {}{^*}\sum _{\chi ^*(d)}\mathcal G_{\chi }(s,w)\sum _{\mathbf l,\mathbf n}a_{\mathbf l,\mathbf n}\chi ^*(-NL)\notag 
\\ &=:&\sum _{d|q}\left (\frac {d}{\pi }\right )^{1-s-w}f_{s,w}^*(d)\label {f}
\end {eqnarray}
and
\begin {eqnarray}\label {g}
\mathcal G_{\chi }(s,w)&=&\frac {\Gamma \left (\frac {1-s+k}{2}\right )}{\Gamma \left (\frac {s+k}{2}\right )}\frac {\Gamma \left (\frac {1-w+k}{2}\right )}{\Gamma \left (\frac {w+k}{2}\right )}=:\left \{ \begin {array}{ll}E_{s,w}&\chi^* \text { even }\\ O_{s,w}&\chi ^*\text { odd }\end {array}\right \} .
\end {eqnarray}
For $(nm,d)=1$ it is straightfoward to establish 
\begin {eqnarray*}
2\sideset {}{^*}\sum _{\chi (d)\atop {\chi ^*(-1)=(-1)^a}}\chi (n\overline m)&=&\sum _{k|d\atop {n\equiv m(k)}}\phi (k)\mu (d/k)+(-1)^a\sum _{k|d\atop {n\equiv -m(k)}}\phi (k)\mu (d/k)
\end {eqnarray*}
so 
\begin {eqnarray}\label {mwyarduon}
f_{s,w}^*(q)&=&E_{s,w}\sum _{\mathbf l,\mathbf n}a_{\mathbf l,\mathbf n}\sideset {}{^*}\sum _{\chi \atop {\text {even}}}\chi ^*(NL)-O_{s,w}\sum _{\mathbf l,\mathbf n}a_{\mathbf l,\mathbf n}\sideset {}{^*}\sum _{\chi \atop {\text {odd}}}\chi ^*(NL)\hspace {10mm}(\text {from \eqref {f}, \eqref {g}})\notag 
\\ &=&E_{s,w}\sum _{k|d}\phi (k)\mu (d/k)\left (\sum _{\mathbf l,\mathbf n\atop {n_1n_4n_6l_3l_4\equiv n_2n_3n_5l_1l_2(k)}}+\sum _{\mathbf l,\mathbf n\atop {n_1n_4n_6l_3l_4\equiv -n_2n_3n_5l_1l_2(k)}}\right )a_{\mathbf l,\mathbf n}\hspace {5mm}(\text {from \eqref {f}, \eqref {g}})\notag 
\\ &&-\hspace {4mm}O_{s,w}\sum _{k|d}\phi (k)\mu (d/k)\left (\sum _{\mathbf l,\mathbf n\atop {n_1n_4n_6l_3l_4\equiv n_2n_3n_5l_1l_2(k)}}-\sum _{\mathbf l,\mathbf n\atop {n_1n_4n_6l_3l_4\equiv -n_2n_3n_5l_1l_2(k)}}\right )a_{\mathbf l,\mathbf n}\notag 
\\ &=&\frac {E_{s,w}}{2}\sum _{k|d}\phi (k)\mu (d/k)\sum _{N,M=1\atop {N\equiv M(k)\atop {(NM,d)=1}}}^\infty \frac {1}{N^{w+1}M^{s+1}}
\notag 
\\ &&\times \hspace {0mm}\sum _{\mathbf l,\mathbf n\atop {\mathbf n\leq d^{1+2\epsilon }\atop {\mathbf l|q\atop {n_1n_4n_6l_3l_4=N\atop {n_2n_3n_5l_1l_2=M}}}}}\mu (\mathbf l)l_1l_3R_\Delta (n_5)R_\Delta (n_6)\chi _0(n_5n_6)(n_1n_3)^{s+w}n_5n_6\mathcal W_{s,w}(\mathbf n)\notag 
\\ &&\hspace {10mm}+\hspace {4mm}\underbrace {\frac {E_{s,w}}{2}\sum _{N\equiv -M}...\hspace {4mm}-\hspace {4mm}\frac {O_{s,w}}{2}\sum _{N\equiv M}...\hspace {4mm}+\hspace {4mm}\frac {O_{s,w}}{2}\sum _{N\equiv -M}...}_{=:\text { similar }}\notag 
\\ &=:&\frac {E_{s,w}}{2}\sum _{k|d}\phi (k)\mu (d/k)\sum _{N,M=1\atop {N\equiv M(k)\atop {(NM,d)=1}}}^\infty \frac {H_q^{s+w}(N,M)}{N^{w+1}M^{s+1}}+\text { similar }
\end {eqnarray}
and we conclude from \eqref {f} 
\begin {eqnarray}\label {rostock}
f_{s,w}(q)&=&\frac {E_{s,w}}{2}\sum _{d|q}\left (\frac {d}{\pi }\right )^{1-s-w}\sum _{k|d}\phi (k)\mu (d/k)\sum _{N,M=1\atop {N\equiv M(k)\atop {(NM,d)=1}}}^\infty \frac {H_q^{s+w}(N,M)}{N^{w+1}M^{s+1}}+\text { similar }\notag 
\\ &=&\frac {E_{s,w}}{2}\sum _{d|q}\left (\frac {d}{\pi }\right )^{1-s-w}\sum _{k|d}\phi (k)\mu (d/k)\left (\sum _{N=M}+\sum _{N\not =M}\right )+\text { similar }\notag 
\\ &=:&\frac {E_{s,w}}{2}\Big (\mathcal D_q(s,w)+\mathcal N_{s,w}(q)\Big )+\text { similar. }
\end {eqnarray}
From \eqref {diagonal} and \eqref {hwyr}
\begin {eqnarray*}
\mathcal W_{s,w}(\mathbf n)&=&
1+\mathcal O\left (d^{100\epsilon }\frac {n_1+n_2+n_3+n_4}{d^{}}\right )\hspace {10mm}\text {for }\mathbf n\leq d^{1+2\epsilon }
\end {eqnarray*}
so from \eqref {mwyarduon} (and note $R_\Delta (n)n\ll 1$ from \eqref {extension} and \eqref {s})
\begin {eqnarray*}
\sum _{N,M=1\atop {(NM,d)=1\atop {N=M}}}^\infty \frac {H_q^{s+w}(N,M)}{N^{w+1}M^{s+1}}&=&\sum _{N,M=1\atop {(NM,d)=1\atop {N=M}}}^\infty \frac {1}{N^{w+1}M^{s+1}}\sum _{\mathbf l,\mathbf n\atop {\mathbf n\leq d^{1+2\epsilon }\atop {\mathbf l|q\atop {n_1n_4n_6l_3l_4=N\atop {n_2n_3n_5l_1l_2=M}}}}}\underbrace {\mu (\mathbf l)l_1l_3R_\Delta (n_5)R_\Delta (n_6)\chi _0(n_5n_6)(n_1n_3)^{s+w}n_5n_6}_{=:\hspace {1mm}b_{\mathbf l,\mathbf n}^{s+w}}\notag 
\\ &&+\hspace {4mm}\mathcal O\left (d^{101\epsilon }\sum _{N=1}^\infty \frac {1}{N}\sum _{l_1,l_3,l_4,n_1,n_4,n_6\atop {\mathbf n\leq d^{1+\epsilon }\atop {n_1n_4n_6l_3l_4=N\atop {l_1|N}}}}\frac {l_1}{n_1n_4n_6}\cdot \frac {n_1+n_4}{d}+...\sum _{...}\frac {l_3}{n_2n_3n_5}\cdot \frac {n_2+n_3}{d}\right )\notag 
\\ &=&\sum _{N=1\atop {(N,d)=1}}^\infty \frac {1}{N^{2+s+w}}\sum _{\mathbf l,\mathbf n\atop {\mathbf l|q\atop {n_1n_4n_6l_3l_4=N\atop {n_2n_3n_5l_1l_2=N}}}}b_{\mathbf l,\mathbf n}^{s+w}+\mathcal O\left (d^{110\epsilon -1}\right )\notag 
\\ &=:&\sum _{N=1\atop {(N,d)=1}}^\infty \frac {h_q^{s+w}(N)^2}{N^{2+s+w}}+\mathcal O\left (d^{110\epsilon -1}\right )
\end {eqnarray*}
so from \eqref {rostock} 
\begin {eqnarray}\label {cinio}
\mathcal D_q(s,w)&=&\sum _{d|q}\left (\frac {d}{\pi }\right )^{1-s-w}\sum _{k|d}\phi (k)\mu (d/k)\sum _{N=1\atop {(N,d)=1}}^\infty \frac {h_q^{s+w}(N)^2}{N^{2+s+w}}+\mathcal O(q^{1+110\epsilon })
\end {eqnarray}
and it is perhaps useful to write down
\begin {eqnarray}\label {hh}
h_q^u(N)&:=&\sum _{abcll'=N\atop {l,l'|q}}\mu (l)\mu (l')lR_\Delta (c)\chi _0(c)a^uc.
\end {eqnarray}
Referring to \eqref {hh} and \eqref {r} straightforward calculations give, for $N$ a power of a prime $p$, 
\begin {eqnarray}\label {h}
h_q^u(N)&=&h_p^u(N)\hspace {1mm}=\hspace {1mm}N^u(1-p^{1-u})\hspace {10mm}\text { if }p|q\notag 
\\ h_q^u(N)&=&h_1^u(N)\hspace {1mm}=\hspace {1mm}\sum _{ac|N}R_\Delta (c)ca^u\hspace {12mm}\text { if }p\nmid q
\end {eqnarray}
so that, where the sum extends over the non-negative powers of $p$,
\begin {eqnarray}\label {hhh}
U:=p^{-u}\hspace {10mm}\mathfrak h_p&:=&\sum _{N}\frac {h_p^u(N)^2}{N^{2+u}}\hspace {1mm}=\hspace {1mm}\frac {p^2U(1-\frac {2}{p}+U)}{p^2U-1}\notag 
\\ \mathfrak h_1&:=&\sum _N\frac {h_1^u(N)^2}{N^{2+u}}\hspace {1mm}=\hspace {1mm}\frac {p^2U(1+\frac {2r}{p}+r^2U)}{p^2U-1}\hspace  {10mm}\text {if }p\nmid \Delta \notag 
\\ \mathfrak h_1^*&:=&\sum _N\frac {h_1^u(N)^2}{N^{2+u}}\hspace {1mm}=\hspace {1mm}1\hspace {10mm}\text {if }p|\Delta 
\end {eqnarray}
(the $\mathfrak h_1$ calculation being admittedly quite lengthy) and therefore
\begin {eqnarray*}\label {glaw}
\sum _{N=1\atop {(N,d)=1}}^\infty \frac {h_q^u(N)^2}{N^{2+u}}&=&
\prod _{p\nmid \Delta }\mathfrak h_1\prod _{p|\Delta }\mathfrak h_1^*\prod _{p|q}\frac {\mathfrak h_p}{\mathfrak h_1}\prod _{p|d}\frac {1}{\mathfrak h_p}\notag 
\\ &=&\zeta (2-s)\prod _{p\nmid \Delta }\left (1+\frac {2r}{p}+r^2U\right )\prod _{p|q}\frac {1-\frac {2}{p}+U}{1+\frac {2r}{p}+r^2U}\prod _{p|d}\frac {p^2U-1}{p^2U\left (1-\frac {2}{p}+U\right )}\notag 
\\ &=:&\zeta (2-u)\mathcal P_\Delta (u-1)\psi _1^u(q)\psi _{2}^u(d)\hspace {10mm}(\text {with $\mathcal P$ as in \eqref {p}})
\end {eqnarray*}
so from \eqref {cinio}, assuming $q$ is square-free,
\begin {eqnarray}\label {dd}
\mathcal D_q(s,w)&=&\zeta (2-s-w)\mathcal P_\Delta (s+w-1)\psi _1^{s+w-1}(q)\sum _{d|q}\left (\frac {d}{\pi }\right )^{1-s-w}\psi _2^{s+w-1}(d)\sum _{k|d}\phi (k)\mu (d/k)+\mathcal O(q^{1+110\epsilon })\notag 
\\ &=&\underbrace {\frac {\zeta (2-s-w)\mathcal P_\Delta (s+w)}{\pi ^{1-s-w}}}_{=:\hspace {1mm}\mathcal Z(s+w)}\frac {\phi (q)\theta _{s+w-1}(q)}{q^{s+w-1}}+\mathcal O(q^{1+110\epsilon })
\end {eqnarray}
where $\theta _s(q)$ is as in \eqref {p}. From \eqref {rostock} and \eqref {mwyarduon} (and $\mathcal W_{s,w}(\mathbf n)\ll 1$ from \eqref {diagonal})
\begin {eqnarray*}
\sum _{q\leq X}\frac {g(q)}{\phi (q)}\mathcal N_{s,w}(d)
&\ll &\sum _{q\leq X}\frac {|g(q)|}{\phi (q)}\sum _{d|q}d\sum _{k|d}\phi (k)\sum _{N,M\leq q^{\mathcal O(1)}\atop {N\equiv M(k)\atop {N\not =M}}}\frac {1}{NM}\sum _{l_3|N\atop {l_1|M}}l_1l_3
\hspace {1mm}\ll \hspace {1mm}X^{\epsilon }
\end {eqnarray*}
so from \eqref {rostock}, \eqref {dd} and \eqref {mwyarduon} (picking up a diagonal term only in one of the ``similar" terms in \eqref {mwyarduon})
\begin {eqnarray*}
\sum _{q\leq X}\frac {g(q)}{\phi (q)}f_{s,w}(q)&=&\mathcal Z(s+w)\cdot \frac {E_{s,w}-O_{s,w}}{2}\cdot \sum _{q\leq X}g_\Delta (q)\frac {\theta _{s+w-1}(q)}{q^{s+w-1}}+\mathcal O\left (X^{111\epsilon }\Big (|E_{s,w}|+|O_{s,w}|\Big )\right )\notag 
\\ &=&\mathcal Z(s+w)\cdot \frac {E_{s,w}-O_{s,w}}{2}\cdot R_{s+w-1}+\mathcal O\left (X^{111\epsilon }\Big (|E_{s,w}|+|O_{s,w}|\Big )\right )
\end {eqnarray*}
with $R_s$ as in \eqref {r}, and we conclude from \eqref {dwr}
\begin {eqnarray}\label {eira}
\int \! \! \! \int \mathcal D&=&\frac {1}{2}\int _{(\epsilon )}\int _{(\epsilon )}\frac {\mathcal Z(s+w)E_{s,w}R_{s+w}\Gamma (s)\Gamma (w+1)X^{s+w+3}dwds}{\Gamma (s+w+4)}-...\int \! \! \! \int ...O_{s,w}...\notag 
\\ &&+\hspace {4mm}\mathcal O\left (X^{3+113\epsilon }\int _{(\epsilon )}\int _{(\epsilon )}\left |\frac {\Gamma (s)\Gamma (w+1)}{\Gamma (s+w+4)}\right |\Big (|E_{s,w}|+|O_{s,w}|\Big )dwds\right )\notag 
\\ &=:&\frac {1}{2}\Big (\mathcal I_E-\mathcal I_O\Big )+\mathcal O\left (X^{3+113\epsilon }\right )
\end {eqnarray}
where from \eqref {g} and the substitution $s+w-1=u$ 
\begin {eqnarray}\label {wy}
\mathcal I_E&=&
\int _{(2\epsilon -1)}\frac {\mathcal Z(u+1)R_{u}X^{u+4}}{\Gamma (u+5)}\left (\underbrace {\int _{(\epsilon )}\frac {\Gamma \left (\frac {1-s}{2}\right )\Gamma \left (\frac {s-u}{2}\right )}{\Gamma \left (\frac {s}{2}\right )\Gamma \left (\frac {u+1-s}{2}\right )}\Gamma (s)\Gamma (u+2-s)ds}_{=:\hspace {1mm}j(u)}\right )du\notag 
\\ \mathcal I_O&=&\int _{(2\epsilon -1)}\frac {\mathcal Z(u+1)R_{u}X^{u+4}}{\Gamma (u+5)}\left (\int _{(\epsilon )}\frac {\Gamma \left (1-\frac {s}{2}\right )\Gamma \left (\frac {s-u+1}{2}\right )}{\Gamma \left (\frac {s+1}{2}\right )\Gamma \left (1-\frac {s-u}{2}\right )}\Gamma (s)\Gamma (u+2-s)ds\right )du.
\end {eqnarray}
For $a,a',b,b'\in \mathbb C$ with 
\[ \underbrace {\mathfrak R\mathfrak e(a-1),\mathfrak R\mathfrak (a'-1)}_{=:A}<\underbrace {\mathfrak R\mathfrak e(b),\mathfrak R\mathfrak e(b')}_{=:B}\]
define\footnote {see https://functions.wolfram.com/HypergeometricFunctions/MeijerG/02/} a special case of Meijer's $G$-function by
\[ G_{2,2}^{2,2}\left (\begin {array}{ll}a&a'\\ b&b'\end {array}\right ):=\frac {1}{2\pi i}\int _{(c)}\Gamma (b-s)\Gamma (b'-s)\Gamma (1-a+s)\Gamma (1-a'+s)ds\hspace {10mm}A<c<B\]
(to get our bearings it may be helpful to note that the ``$a,a'$-poles" are to the left of $A$ and the ``$b,b'$-poles" are to the right of $B$) which is known\footnote {see https://functions.wolfram.com/HypergeometricFunctions/MeijerG/03/01/04/21/ and note that Gauss' hypergeometric function ${}_2F_1(a,b,c;z)$ has ${}_2F_1(a,b,c;0)=1$} to satisfy
\begin {eqnarray*}
G_{2,2}^{2,2}\left (\begin {array}{ll}a&a'\\ b&b'\end {array}\right )&=&\frac {\Gamma (b+1-a)\Gamma (b'+1-a)\Gamma (b+1-a')\Gamma (b'+1-a')}{
\Gamma (2+b+b'-a-a')}. 
\end {eqnarray*}
In the following we use freely the formulas $s\Gamma (s)=\Gamma (s+1)$ and $\Gamma (s/2)\Gamma (s/2+1/2)=\sqrt \pi 2^{1-s}\Gamma (s)$. Since the integrand of $j(u)$ is $...$ times
\begin {eqnarray*}
2\Gamma \left (\frac {1}{2}-\frac {s}{2}\right )\Gamma \left (2+\frac {u}{2}-\frac {s}{2}\right )\Gamma \left (\frac {1}{2}+\frac {s}{2}\right )\Gamma \left (-\frac {u}{2}+\frac {s}{2}\right )-\Gamma \left (\frac {1}{2}-\frac {s}{2}\right )\Gamma \left (1+\frac {u}{2}-\frac {s}{2}\right )\Gamma \left (\frac {1}{2}+\frac {s}{2}\right )\Gamma \left (-\frac {u}{2}+\frac {s}{2}\right )
\end {eqnarray*}
we have (with $c=1/4$)
\begin {eqnarray*}
j(u)&=&\frac {2^{s}}{\sqrt \pi }\left \{ 2G_{2,2}^{2,2}\left (\begin {array}{ll}1/2&1+u/2\\ 1/2&2+u/2\end {array}\right )-G_{2,2}^{2,2}\left (\begin {array}{ll}1/2&1+u/2\\ 1/2&1+u/2\end {array}\right )\right \} 
\\ &=&\frac {2^{s}}{\sqrt \pi }\left \{ \Gamma \left (\frac {5}{2}+\frac {u}{2}\right )\Gamma \left (\frac {1}{2}-\frac {u}{2}\right )-\Gamma \left (\frac {3}{2}+\frac {u}{2}\right )\Gamma \left (\frac {1}{2}-\frac {u}{2}\right )\right \} 
\\ &=:&g_E(u)
\end {eqnarray*}
so after a similar calculation with
\begin {eqnarray*}
g_O(u)&:=&\frac {2^{s}}{\sqrt \pi }\Gamma \left (\frac {3}{2}+\frac {u}{2}\right )\Gamma \left (\frac {3}{2}-\frac {u}{2}\right ).
\end {eqnarray*}
we get from \eqref {wy}
\begin {eqnarray}\label {haul}
\mathcal I_E-\mathcal I_O&=&\int _{(2\epsilon -1)}\frac {\mathcal Z(u+1)R_{u}X^{u+4}}{\Gamma (u+5)}\Big (g_E(u)-g_O(u)\Big )du\notag 
\\ &=&\int _{(2\epsilon -1)}\pi ^{s}\frac {\mathcal Z(u+1)R_{u}\Gamma \left (\frac {1}{2}-\frac {u}{2}\right )X^{u+4}}{\Gamma \left (u/2\right )(u+2)(u+3)(u+4)}du
\end {eqnarray}
and now \eqref {1} follows from \eqref {dwr}, \eqref {qq}, \eqref {eira}, and the definitions in \eqref {dd}, \eqref {r}.
\\ 
\\ 
\begin {center}
{\section {- \hspace {3mm}Proof of \eqref {2}}}
\end {center}
All vectors throughout are two-dimensional vectors of natural numbers and we write
\begin {eqnarray*}
\mathbf d|n\hspace {10mm}&&\text {for}\hspace {10mm}d,d'|n
\\ (\mathbf d,n)=1\hspace {10mm}&&\text {for}\hspace {10mm}(dd',n)=1
\\ F(\mathbf d)\hspace {10mm}&&\text {for}\hspace {10mm}F(d)F(d').
\end {eqnarray*}
Recall the definitions of $f_\Delta ,g_\Delta $ from \eqref {s}. As 
\[ \sum _{l+l'=n\atop {\mathbf D|\mathbf l}}1\hspace {1mm}=\hspace {1mm}\underbrace {\frac {n}{[D,D']}}_{(D,D')|n}\hspace {1mm}+\hspace {1mm}\mathcal O(1)\hspace {10mm}\text { and from \eqref {bound} }\hspace {10mm}\sum _{[D,D]>X}\frac {g(\mathbf D)}{[D,D']}\hspace {1mm}\ll \hspace {1mm}X^{\epsilon -1}\]
we have, picking up the coprimality conditions with the M\" obius function,
\begin {eqnarray}
\sum _{l+l'\leq X\atop {(\mathbf l,q)=1\atop {\tilde l\equiv 0(q)}}}(X-\tilde l)^2\tilde lf_\Delta (\mathbf l)&=&
\sum _{n\leq X\atop {n\equiv 0(q)}}(X-n)^2n^2\underbrace {\sum _{\mathbf h|q}\mu (\mathbf h)\sum _{\mathbf D\atop {([D,h],[D',h'])|n}}\frac {g_\Delta (\mathbf D)}{[D,h,D',h']}}_{=:a_q(n)}+\mathcal O\left (\frac {X^{4+\epsilon }}{q}\right )\notag 
\\ &=&2\int _{(1)}\frac {\mathcal K_{q}(s)X^{s+4}ds}{(s+2)(s+3)(s+4)}+\mathcal O\left (\frac {X^{4+\epsilon }}{q}\right )\label {k}
\\ \text {where }\hspace {10mm}\mathcal K_q(s)&:=&\sum _{n=1\atop {n\equiv 0(q)}}^\infty \frac {a_q(n)}{n^s}\label {kk}
\end {eqnarray}
with a standard Perron formula. Write $(f_1\star _Lf_2)(n):=\sum _{[a,b]=n}f_1(a)f_2(b)$. Then 
\begin {eqnarray*}
m_q(A)\hspace {1mm}:=\hspace {1mm}\sum _{[D,h]=A\atop {h|q}}\mu (h)g_\Delta (D)\hspace {1mm}=\hspace {1mm}\prod _{p^\alpha ||A\atop {p\nmid q}}g_\Delta (p^\alpha )\prod _{p^\alpha ||A\atop {p|q}}(g_\Delta\star _L\mu )(p^\alpha )\hspace {1mm}=\hspace {1mm}g_\Delta (A)\underbrace {\prod _{p^\alpha ||A\atop {p|q}}\frac {(g_\Delta \star _L\mu )(p^\alpha )}{g_\Delta (p^\alpha )}}_{=:G_q(A)}
\end {eqnarray*}
so, after a straightforward calculation,
\begin {eqnarray*}
a_q(n)&=&\sum _{A,A'=1\atop {(A,A')|n}}^\infty \frac {m_q(\mathbf A)}{[A,A']}
\\ &=&\prod _{p\nmid n}\left (1+\frac {2g_\Delta (p)G_q(p)}{p}\right )\prod _{p|n}\left (1+\frac {2g_\Delta (p)G_q(p)}{p}+\frac {g_\Delta (p)^2G_q(p)^2}{p}\right )
\\ &=:&\prod _{p\nmid n}u_q(p)\prod _{p|n}\left (u_q(p)+U_q(p)\right )
\\ &=&\prod _pu_1(p)\prod _{p|q}\frac {u_p(p)+U_p(p)}{u_1(p)}\prod _{p|n\atop {p\nmid q}}\left (1+\frac {U_1(p)}{u_1(p)}\right )
\\ &=:&CK(q)J_q(n)
\end {eqnarray*}
so that from \eqref {kk}
\begin {eqnarray}\label {papur}
\mathcal K_q(s)&=&\frac {CK(q)}{q^s}\sum _{n=1}^\infty \frac {J_q(n)}{n^s}\notag 
\\ &=&\frac {C\zeta (s)K(q)}{q^s}\prod _{p}\Big (1+Y\left (J_1(p)-1\right )\Big )\prod _{p|q}\frac {1}{1+Y\left (J_1(p)-1\right )}\hspace {10mm}Y:=p^{-s}\notag 
\\ &=&\frac {\zeta (s)}{q^s}\prod _{p}\Big (u_1(p)+YU_1(p)\Big )\prod _{p|q}\frac {u_p(p)+U_p(p)}{u_1(p)+YU_1(p)}\notag 
\\ &=&\frac {\zeta (s)}{q^s}\prod _{p\nmid \Delta }\left (1+\frac {2g_1(p)}{p}+\frac {Yg_1(p)^2}{p}\right )\prod _{p|q}\frac {1-1/p}{1+\frac {2g_1(p)}{p}+\frac {Yg_1(p)^2}{p}}\hspace {10mm}\text {(from \eqref {bound} and $g_1(p)G_p(p)=-1$)}\notag \notag 
\\ &=&\zeta (s)\mathcal P_\Delta (s)\frac {\theta _{s}(q)}{q^s}\hspace {10mm}\text {(from \eqref {bound} and \eqref {p})}\notag 
\\ &=:&\zeta (s)\mathcal P_\Delta (s)F_s(q).
\end {eqnarray}
Putting this in \eqref {k} and recalling \eqref {ll} gives 
\begin {eqnarray*}
\mathcal L_\Delta (q)&=&
\mathcal P_\Delta (1)F_1(q)X^5+2\int _{(1/2)}\frac {\zeta (s)F_s(q)X^{s+4}ds}{(s+2)(s+3)(s+4)}+\mathcal O\left (\frac {X^{4+\epsilon }}{q}\right )
\end {eqnarray*}
so
\begin {eqnarray*}
\sum _{q>X}g_\Delta (q)\mathcal L_\Delta (q)&=&X^5\sum _{q>X}g_\Delta (q)F_1(q)+2\int _{(1/2)}\overbrace {\zeta (s)\mathcal P_\Delta (s)\left (\sum _{q=1}^\infty g_\Delta (q)F_s(q)\right )}^{(\star )}\frac {X^{s+4}ds}{(s+2)(s+3)(s+4)}
\\ &&-\hspace {4mm}2\int _{(1/2)}\zeta (s)\mathcal P_\Delta (s)\left (\sum _{q\leq X}g_\Delta (q)F_{q}(s)\right )\frac {X^{s+4}ds}{(s+2)(s+3)(s+4)}+\mathcal O\left (X^{4+\epsilon }\sum _{q>X}\frac {|g_\Delta (q)|}{q}\right ).
\end {eqnarray*}
Recall the definitions of $\mathcal A^>(X),\mathcal Q_\Delta ^>(s),\mathcal E_\Delta ^>$ from \eqref {a}, \eqref {q}, \eqref {e}. From \eqref {papur} the term $(\star )$ is $\mathcal Q_\Delta ^>(s)$ so the first, second and fourth terms together are 
\[ X^5\mathcal A^>(X)+\mathcal E_\Delta ^>+\mathcal O\left (X^{3+\epsilon }\right ).\]
Recall $R_s,\mathcal R_\Delta $ from \eqref {r}. From \eqref {papur} the $q$-sum in the last integral is $R_s$ so this integral, after moving it to the left and using the functional equation for the Riemann zeta function, is $\mathcal R_\Delta $. So we get \eqref {2}.


\begin {center}

\end {center}
\hspace {1mm}
\\
\\
\\
\\
\\  
\\ \emph {Tomos Parry
\\ Bilkent University, Ankara, Turkey
\\ tomos.parry1729@hotmail.co.uk}
\\
\\\ There is no conflict of interest in my submission and there is no associated data with the submission.

\end {document}